\documentclass[a4paper,10pt]{elsarticle}
\usepackage{framed,multirow}
\usepackage[utf8]{inputenc}
\usepackage{bm}
\usepackage{amsmath}
\usepackage{mathrsfs}
\usepackage{graphicx}
\usepackage{epsfig, setspace}
\usepackage{url}
\usepackage{graphicx, transparent, color}
\usepackage{algorithm} 
\usepackage{algorithmicx}
\usepackage{etex}
\usepackage[bookmarks=false]{hyperref}
\usepackage{amsmath}
\usepackage{xcolor}
\usepackage{rotating} 
\usepackage{pdflscape} 
\usepackage{silence}
\usepackage{ulem,cancel}
\usepackage[utf8]{inputenc}

\usepackage{graphicx}
\usepackage{float}
\usepackage{subfigure}
\usepackage{subfigmat}

\usepackage{tikz}
\usepackage{capt-of}
\usepackage{setspace,lipsum}
\usepackage{lineno}

\usepackage[top=18truemm,bottom=18truemm,left=20truemm,right=20truemm]{geometry}

\begin{document}
\begin{frontmatter}
\title{\textcolor{black}{On the Importance of High-Frequency Damping in High-Order Conservative Finite-Difference Schemes for Viscous Fluxes}}
\author[AA_address]{Amareshwara Sainadh Chamarthi \cortext[cor1]{Corresponding author. \\ 
E-mail address: skywayman9@gmail.com (Amareshwara Sainadh  Ch.).}}
\author[AA_address]{Sean Bokor} 
\author[AA_address]{Steven H.\ Frankel}
\address[AA_address]{Faculty of Mechanical Engineering, Technion - Israel Institute of Technology, Haifa, Israel}

\begin{abstract}

\textcolor{black}{
This paper discusses the importance of high-frequency damping in high-order conservative finite-difference schemes for viscous terms in the Navier-Stokes equations. Investigating nonlinear instability  encountered in a high-resolution viscous shock-tube simulation, we have discovered that a modification to the viscous scheme rather than the inviscid scheme resolves a problem with spurious oscillations around shocks. The modification introduces a term responsible for high-frequency damping that is missing in a conservative high-order viscous scheme. The importance of damping has been known for schemes designed for unstructured grids. However, it has not been recognized well in very high-order difference schemes, especially in conservative difference schemes. Here, we discuss how it is easily missed in a conservative scheme and how to improve such schemes by a suitably designed damping term. 
}

\end{abstract}


\end{frontmatter}

\section{Introduction}\label{sec-1}
\textcolor{black}{
Historically, the discretization of the viscous terms in the Navier-Stokes equations have received less attention in the literature than the discretization of the inviscid terms. Much of the research is focused on the discretization of the inviscid terms especially to reduce spurious oscillations near discontinuities. While non-oscillatory inviscid schemes are essential, we have noted that viscous schemes also play a key role in the quality of the simulation, especially regarding the prevention of spurious oscillations near discontinuities in high-resolution viscous shock-tube simulations. It is found particularly important for high-order conservative difference schemes, e.g., in a general form in the $x$-direction, 
\begin{eqnarray}
 \left.
  \frac{\partial F}{\partial x} \right|_{j} 
  = \frac{1}{\Delta x}
   \left[   
   C_{j-1} \left(   \hat{F}_{j+1/2}   -  \hat{F}_{j-1/2}  \right) 
+  C_{j}  \left(   \hat{F}_{j+3/2}   -  \hat{F}_{j-3/2}    \right) 
+  C_{j+1}  \left(   \hat{F}_{j+5/2}   -  \hat{F}_{j-5/2}    \right) 
    \right],
    \label{bad_sixth-order_scheme}
\end{eqnarray}
which approximates the flux divergence $\partial F / \partial x$ at the center of a cell $j$, where $\hat{F}$ denotes a numerical flux and point-valued solutions are stored at the cell centers, $u_j$, $ j = 1,2,3, \cdots$. For a diffusion term (a representative of the viscous terms), $\partial F/  \partial x = \partial^2 u / \partial x^2 $ with $F =  \partial u / \partial x$, the numerical flux needs to approximate the solution gradient $\partial u / \partial x$, and is typically designed for a desired order of accuracy. For example, Shen et al. \cite{shen2010large} developed a sixth-order scheme by setting $C_{j-1} = 75/64$, $C_{j} = -25/384$, $C_{j+1} = 3/640$, and defining the numerical flux as 
 \begin{eqnarray}
 \hat{F}_{j+1/2}    = 
    \frac{3}{256} \left(   \frac{ \partial u }{  \partial x }   \right)_{j-2} 
-  \frac{25}{256} \left(   \frac{ \partial u }{  \partial x }   \right)_{j-1} 
+  \frac{150}{256} \left(   \frac{ \partial u }{  \partial x }   \right)_{j} 
+  \frac{150}{256} \left(   \frac{ \partial u }{  \partial x }   \right)_{j+1} 
-  \frac{25}{256} \left(   \frac{ \partial u }{  \partial x }   \right)_{j+2} 
+  \frac{3}{256} \left(   \frac{ \partial u }{  \partial x }   \right)_{j+3} ,
 \label{bad_6th-order_scheme_flux}
\end{eqnarray}
where each solution gradient is computed by 
\begin{eqnarray} 
 \left(   \frac{ \partial u }{  \partial x }   \right)_{j}  
 =
 \frac{  1 }{60} \left[ 
         45  \frac{  u_{j+1} - u_{j-1}   }{\Delta x}
        -  9  \frac{  u_{j+2} - u_{j-2}   }{\Delta x}
     +   {\color{black}   \frac{  u_{j+3} - u_{j-3}   }{\Delta x} }
 \right].
 \end{eqnarray}
After employing this scheme (referred to as Shen's scheme from here on)}, for the viscous discretization in our high-order Navier-Stokes code, \textcolor{black}{we} encountered {\color{black} spurious oscillations} around a shock in a viscous shock-tube simulation \textcolor{black}{(see Fig \ref{fig:500_shen_nodamp_50})}. Typically, one would attempt to resolve this issue by modifying the inviscid algorithm, but we have discovered that the problem can be resolved by improving the viscous discretization. It will be shown later that the above viscous scheme is missing a high-frequency damping property (or, in other words, is not $h$-elliptic \cite{trottenberg2000multigrid}) which allows the odd-even decoupling error modes to persist, which can cause serious problems especially around shock waves, where high-frequency modes dominate. As we will show, introducing high-frequency damping resolves this particular problem. However, rather than just adding a damping term to the above scheme, we propose a simple construction of a high-order conservative scheme for viscous terms, which explicitly introduces the damping property and results in a sixth-order scheme with a significantly smaller stencil. 

This short note's objective is to discuss an approach to constructing a class of conservative viscous schemes equipped with a high-frequency damping mechanism suitable for simulations with shock waves. It will be shown that the proposed scheme provides the damping mechanism and significantly improves the resolution of the numerical results itself. The rest of the note is organized as follows. Section \ref{sec-2} introduces the governing equations and the discretization of viscous fluxes is presented in Section \ref{sec-3}. Numerical results are presented in Section \ref{sec-4}, and finally, in Section \ref{sec-5}, we provide conclusions.


\section{\color{black} Governing Equations and Target Conservative Finite-Difference Scheme}\label{sec-2}
The compressible Navier-Stokes (NS) equations governing single-component fluid flows in a two-dimensional (2-D) Cartesian coordinate system without loss of generality can be expressed as:

\begin{equation}\label{CNS-base}
\frac{\partial \mathbf{Q}}{\partial t}+\frac{\partial \mathbf{F^c}}{\partial x}+\frac{\partial \mathbf{G^c}}{\partial y}+\frac{\partial \mathbf{F^v}}{\partial x}+\frac{\partial \mathbf{G^v}}{\partial y}= {\color{black} \mathbf{0}},
\end{equation}
where  $\textbf{Q}$ is the conservative variable vector, $\mathbf{F^c}$, $\mathbf{G^c}$, and $\mathbf{F^v}$, $\mathbf{G^v}$, are the convective (superscript $c$) and viscous (superscript $v$) flux vectors in each coordinate direction, respectively. The conservative variable, convective, and viscous flux vectors are given as:
\begin{equation}
\mathbf{Q}=\left[\begin{array}{c}
\rho \\
\rho u \\
\rho v \\
E
\end{array}\right],\;
\mathbf{F^c}=\left[\begin{array}{c}
\rho u \\
\rho u^{2}+p \\
\rho u v \\
(E+p) u
\end{array}\right],\;
\mathbf{G^c}=\left[\begin{array}{c}
\rho v \\
\rho u v \\
\rho v^{2}+p \\
(E+p) v
\end{array}\right], 
\mathbf{F^v}=\left[\begin{array}{c}
0 \\
-\tau_{x x} \\
-\tau_{y x} \\
- {\color{black}  \tau_{xu} } +q_{x}
\end{array}\right],\; 
\mathbf{G^v}=\left[\begin{array}{c}
0 \\
-\tau_{x y} \\
-\tau_{y y} \\
-  {\color{black}  \tau_{yu} } +q_{y}
\end{array}\right],
\end{equation}
where $\rho$ is the density, $u$ and $v$ are velocity components in the $x-$ and $y-$ directions, respectively,  {\color{black} $\tau_{xu} = \tau_{x x} u + \tau_{x y} v$, $ \tau_{yu} = \tau_{y x} u + \tau_{y y} v$, $E$ is the total energy per unit volume, $p = (\gamma -1) \left( E - \rho \frac{(u^2+v^2)}{2} \right)$ is the pressure with $\gamma=7/5$, which closes the system.}
{\color{black} The viscous stresses} and the heat flux $q$ are given by:
\begin{equation}\label{eqn:5-stress}
\tau_{x x}=\frac{2}{3} \mu\left(2 \frac{\partial u}{\partial x}-\frac{\partial v}{\partial y}\right), \quad \tau_{x y}=\tau_{y x}=\mu\left(\frac{\partial u}{\partial y}+\frac{\partial v}{\partial x}\right), \quad {\color{black}\tau_{y y}=\frac{2}{3} \mu\left(2 \frac{\partial v}{\partial y}-\frac{\partial u}{\partial x}\right),}
\end{equation}
\begin{equation}\label{eqn:6-heat}
\begin{aligned}
q_{x} &=  - {\color{black} \kappa} \frac{\partial T}{\partial x}, \quad q_{y} =  - {\color{black} \kappa} \frac{\partial T}{\partial y},
{\color{black}  \quad T=\frac{\gamma p}{(\gamma-1) C_{p} \rho}  },
\end{aligned}
\end{equation}
where {\color{black} $T$ is the temperature}, $\mu$ is the dynamic viscosity, {\color{black} $\kappa$ is the heat conductivity given by $\kappa$ = $C_p\mu/Pr$, $Pr$ is the Prandtl number, $C_p$ is the specific heat at constant pressure. For the test case considered in this note,} \textcolor{black}{we assume air as an ideal gas, $Pr=0.73$, and $C_p = R \gamma / (\gamma - 1)$ where $R$ is the gas constant.}
 
We apply a conservative numerical method to solve \textcolor{black}{the given system of} equations. The time evolution of the vector of cell-centered conservative variables $\mathbf{\hat Q}$ is given by the following semi-discrete equation applied to a Cartesian grid cell $I_{j,i} = [x_{j-\frac{1}{2}}, x_{j+\frac{1}{2}}] \times [y_{i-\frac{1}{2}}, y_{i+\frac{1}{2}}]$, expressed as {\color{black} a semi-discrete system}:

\begin{equation}\label{eqn-differencing}
\begin{aligned}
\frac{\mathrm{d}}{\mathrm{dt}} {\mathbf{\hat Q}}_{j,i}
= \mathbf{Res}_{j,i}, 
\end{aligned}
\end{equation}
where {\color{black} $ {\mathbf{\hat Q}}_{j,i}$ denotes the vector of numerical solutions at the cell $(j,i)$,} and $ \mathbf{Res}_{j,i}$ denotes the residual:
\begin{equation}\label{eqn-differencing_residual}
\begin{aligned}
\mathbf{Res}_{j,i}=&
{\color{black} 
-\frac{
\left(\mathbf {\hat{F}^c}_{j+ \frac{1}{2}, i}-\mathbf {\hat{F}^c}_{j- \frac{1}{2}, i}\right)-\left(\mathbf {\hat{F}^v}_{j+ \frac{1}{2}, i}
-\mathbf {\hat{F}^v}_{j- \frac{1}{2}, i}\right)}{\Delta x}
-\frac{ \left(\mathbf {\hat{G}^c}_{i+ \frac{1}{2}, j}-\mathbf {\hat{G}^c}_{i- \frac{1}{2}, j}\right)-\left(\mathbf {\hat{G}^v}_{i+ \frac{1}{2}, j}-\mathbf {\hat{G}^v}_{i- \frac{1}{2}, j}\right)}{\Delta y} ,}
\end{aligned}
\end{equation}
where $\Delta x=x_{j+1 / 2}-x_{j-1 / 2}$, $\Delta y=y_{i+1 / 2}-y_{i-1 / 2}$, and $\mathbf {\hat{F}^c}$, $\mathbf {\hat{G}^c}$ and $\mathbf {\hat{F}^v}$, $\mathbf {\hat{G}^v}$ are interpreted as numerical approximations of the convective and viscous fluxes in the $x-$, and $y-$directions, respectively. The inviscid fluxes are computed using the \textcolor{black}{Monotonicity Preserving (MP)} scheme of Suresh and Hyunh \cite{suresh1997accurate}, where we interpolate characteristic variables projected from the primitive variables, and the component-wise local Lax-Friedrichs (cLLF) Riemann solver is used.  {\color{black}  For the viscous terms, the above form is modified as in Equation (\ref{bad_sixth-order_scheme}) for the sixth-order scheme of Shen's scheme, but a simpler sixth-order scheme can be constructed in the above form as we discuss in the next section. 

The semi-discrete system (\ref{eqn-differencing}) is then integrated in time by the third-order TVD Runge-Kutta scheme \cite{jiang1995}. {\color{black} The maximum time step allowable for explicit time-stepping for viscous flow simulations is determined by:}
\begin{equation}\label{eqn:cfl}
\Delta t= \text{CFL} \cdot \min \left(\Delta t_{viscous}, \Delta t_{inviscid}\right),
\end{equation}
where the CFL number is fixed at $0.2$ for all the simulations presented here, 
\begin{equation}
\Delta t_{inviscid}=  \min _{j, i}\left(\frac{\Delta x_{i}}{\left|u_{j, i}\right|+c_{j, i}}, \frac{\Delta y_{j}}{\left|v_{j, i}\right|+c_{j, i}}\right),
 \quad 
 {\color{black}\Delta t_{viscous}=  \min _{j, i}\left(\frac{1}{\alpha} \frac{\Delta x_{j}^{2}}{\nu_{j, i}},\frac{1}{\alpha}  \frac{\Delta y_{i}^{2}}{\nu_{j, i}}\right)},
\end{equation}
$c$ is the speed of sound and given by $c=\sqrt{\gamma{p}/\rho}$, $\alpha$ is a parameter associated with the damping property of a viscous scheme as discussed in the next section, \textcolor{black}{ and $\nu$ is the kinematic viscosity defined as $\nu = \mu / \rho$.} We have found the viscous time step restriction is important for stable computations. Nishikawa derived the above viscous time step restriction in Ref.\textcolor{black}{\cite{nishikawa:AIAA2010}}; see Equation 4.22 there, which is applicable here as well. This time step restriction is similar but much less stringent than a fourth-order discontinuous Galerkin scheme for diffusion. See Table 6.1 of Ref.\textcolor{black}{\cite{nishikawa:AIAA2010}}. 
 }

\section{Efficient Construction of Sixth-Order Viscous Scheme with {\color{black} High-Frequency Damping}}\label{sec-3}
In this section, we describe the discretization of the viscous fluxes. For simplicity, we will consider only a one-dimensional scenario where the viscous flux at the interface can be computed as:
\begin{equation}
\mathbf{\hat F^v}_{j+\frac{1}{2}}=\left[\begin{array}{c}
0 \\
-\tau_{j+\frac{1}{2}} \\
-\tau_{j+\frac{1}{2}} u_{j+\frac{1}{2}}+q_{j+\frac{1}{2}}
\end{array}\right],
\end{equation}
where
\begin{equation}\label{eqn:visc-interface}
\tau_{j+1 / 2}=\frac{4}{3} \mu \left(\frac{\partial u}{\partial x}\right)_{j+1 / 2}, \quad q_{j+1 / 2}=-\frac{\mu}{\operatorname{Pr}(\gamma-1)}\left(\frac{\partial T}{\partial x}\right)_{j+1 / 2},
\end{equation}
where $u_{j+\frac{1}{2}}$ is evaluated with high-order solutions reconstructed at the face. As it can be seen from Equation (\ref{eqn:visc-interface}), similar to the inviscid fluxes, the viscous fluxes at the cell interfaces, ${j+\frac{1}{2}}$, have to be evaluated. In Ref. \textcolor{black}{\cite{nishikawa:AIAA2010}}, Nishikawa presented a general technique \textcolor{black}{for the design of} diffusion schemes with high-frequency damping properties. \textcolor{black}{The numerical flux for his diffusion scheme} consists of two essential parts: a consistent term and a damping term. The former approximates the diffusion flux consistently and accurately, while the latter has an order property (i.e., vanishes eventually in the grid refinement) and plays the role of high-frequency damping. In a typical difference scheme, the damping term is proportional to the solution jump at an interface, $(u_R - u_L)/\Delta x$, where $u_L$ and $u_R$ are reconstructed solutions from the left and right cells, respectively. {\color{black} To introduce high-frequency damping to the sixth-order scheme (\ref{bad_sixth-order_scheme}), we can simply add the damping term to the numerical flux (\ref{bad_6th-order_scheme_flux}). But to preserve sixth-order accuracy, the reconstructed solutions would have to satisfy $u_R - u_L = O(\Delta x^7)$; it is effective but may not be the most efficient approach. A more economical approach would be to construct a high-order accurate gradient at an interface by combining lower-order consistent and damping terms: 
\begin{equation} \label{good_6tth-order_scheme_flux}
\begin{aligned}
  \left(\frac{\partial u}{\partial x}\right)_{j+1 / 2}=\underbrace{\frac{1}{2}\left[\left(\frac{\partial u}{\partial x}\right)_{j}+\left(\frac{\partial u}{\partial x}\right)_{j+1}\right]}_{\text{Consistent term}}+\underbrace{\frac{\alpha}{ 2 \Delta x}\left(u_R- u_L\right)}_{\text{Damping term}}, \\
\end{aligned}
\end{equation}
where the gradients at the cell-centers are computed with the fourth-order central differences:
\begin{eqnarray} 
 \left(   \frac{ \partial u }{  \partial x }   \right)_{j}  =  \frac{8\left({ {u}}_{j+1}-{ {u}}_{j-1}\right)- \left({ {u}}_{j+2}-{ {u}}_{j-2}\right) }{12 \Delta x},
\quad
 \left(   \frac{ \partial u }{  \partial x }   \right)_{j+1}   =  \frac{8\left({ {u}}_{j+2}-{ {u}}_{j}\right)- \left({ {u}}_{j+3}-{ {u}}_{j-1}\right) }{12 \Delta x},
 \label{fourth_order_central_alpha_damp_sixth_uxj_uxjp1}
\end{eqnarray}
and the solutions are reconstructed with quadratic reconstruction: 
\begin{eqnarray} 
 u_L = u_j + \left(   \frac{ \partial u }{  \partial x }   \right)_{j} \frac{\Delta x}{2} + \beta  \left({ {u}}_{j+1} - 2{ {u}}_{j} + { {u}}_{j-1}\right) , \quad 
 u_R =  u_{j+1}  - \left(   \frac{ \partial u }{  \partial x }   \right)_{j+1} \frac{\Delta x}{2} + \beta  \left({ {u}}_{j+2} - 2{ {u}}_{j+1} + { {u}}_{j}\right).
 \label{uLuR_damp_sixth}
\end{eqnarray}
The parameters $\alpha$ and $\beta$ can be determined to achieve sixth-order accuracy in the resulting viscous scheme, as we will discuss below. This approach can be considered as an extension of the construction described in Ref.\textcolor{black}{\cite{nishikawa:AIAA2010}} for deriving fourth-order schemes: second-order central differences are extended to the fourth-order formulas Equation (\ref{fourth_order_central_alpha_damp_sixth_uxj_uxjp1}); and linear reconstructions are extended to quadratic reconstructions Equation (\ref{uLuR_damp_sixth}). Note that this sixth-order scheme has a stencil with seven cells (three on each side of the target cell $j$) whereas the Shen's scheme based on the interface gradients Equation (\ref{bad_6th-order_scheme_flux}) has a stencil of 17 cells. From here on, the above scheme is referred to as the new sixth-order $\alpha$-damping scheme. 
}

\begin{figure}[H]
\centering
 \includegraphics[width=0.46\textwidth]{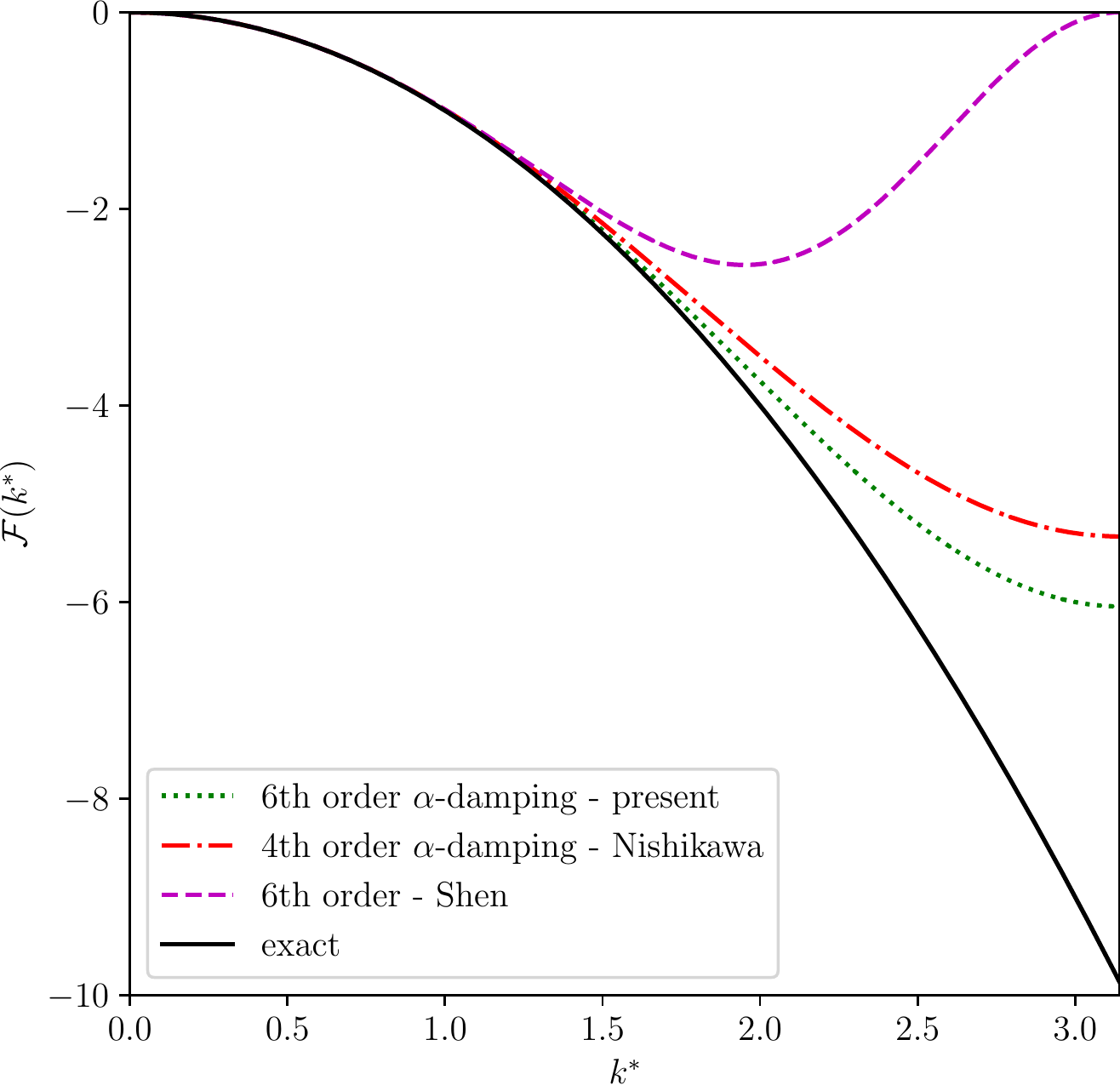}
\caption{Fourier operator for diffusion versus wavenumber.
Differencing error for second derivative vs wavenumber.}
\label{fig:second_deriv}
\end{figure}

\textcolor{black}
{\color{black} A Fourier analysis is presented here to demonstrate sixth-order accuracy and a high-frequency damping property of the proposed scheme. Consider a Fourier mode with a wave number $k$ defined on a grid with $N$ points over a domain $x \in [0,L]$, $k = \frac{2 \pi}{L} n$ is the wave number, and $n = 1, 2, ..., N/2$, 
\begin{equation}
  u(x) = e^{i k x},
\end{equation}
where $i = \sqrt{-1}$, which gives, when applied to the exact diffusion operator, 
\begin{equation}
   \frac{ \partial^2 u }{\partial x^2} =  (i k)^2 e^{i k x} = - k^2 e^{i k x} =   \frac{ -  (k^*)^2 }{\Delta x^2}  e^{i k x}, 
\end{equation} 
showing that the exact diffusion operator is characterized by the negative of the square of the normalized wave number, $- (k^*)^2$, which is plotted as a reference in Fig. \ref{fig:second_deriv}. This can be extended to a discrete form by plugging in values of $x_j$ such that $u_j =   e^{i k x_j}$ with $x_j = \frac{L}{N} j$ 
for $j = 0, 1, 2, ..., N-1$ grid points. Substituting the Fourier mode into a discretization residual, we obtain
\begin{equation}
\left.  \frac{\partial F}{\partial x} \right|_{j}  =  \frac{ {\cal F}(k^*)  }{\Delta x^2}   e^{i k x_j} , 
\end{equation}
where $\left.  \frac{\partial F}{\partial x} \right|_{j}  $ may be given by Equation (\ref{bad_sixth-order_scheme}) or by $ \left.
  \frac{\partial F}{\partial x} \right|_{j} =\left(   \hat{F}_{j+1/2}   -  \hat{F}_{j-1/2}  \right) / \Delta x$ for the new sixth-order scheme, ${\cal F}(k^*)$  is an operator approximating the exact operator $-  (k^*)^2$. For the new sixth-order accurate scheme based on the face gradient (\ref{good_6tth-order_scheme_flux}), we find the operator ${\cal F}(k^*)$ to be}
  
{\color{black} 
\begin{equation}\label{mod-newsixth}
{\cal F}(k^*) =-\frac{1}{6} \sin ^2\left(\frac{k^*}{2}\right) \left[ 24 \alpha  \beta +5 \alpha +6 \left\{ \alpha  (4 \beta -1)+2 \right\} \cos (k^*)+(\alpha -2) \cos (2 k^*)+14 \right].
\end{equation}
By expanding the above equation, we obtain

\begin{equation}
{\cal F}(k^*) = 
- (k^*)^2 \left[    1  - \frac{\alpha (12 \beta -1)+4 }{24}   (k^*)^2 - \frac{9-5 \alpha  (6 \beta +1)}{360}   (k^*)^4 -   \frac{63\alpha(4 \beta +3)-376}{40320} (k^*)^6 + O(  (k^*)^8 )    \right],
\end{equation}

from which we find that the second- and fourth-order errors can be eliminated by $\alpha=38/15$ and $\beta=-11/228$ and thus the scheme achieves sixth-order 
accuracy. Next, we consider Shen's scheme \cite{shen2010large}, Equation (\ref{bad_6th-order_scheme_flux}). Substituting the Fourier mode, we obtain the following equation for the Fourier operator: 
\begin{equation}\label{mod-shen}
  {\cal F}^{Shen}(k^*) = 
  - \frac{   \sin^2(k^*)  \left[     3   \cos^2 (k^*) -14  \cos (k^*) +  43    \right]  
                                  \left[     9   \cos^2 (k^*) -58  \cos (k^*) +  529   \right] 
                                 \left[     2   \cos^2 (k^*) - 9  \cos (k^*) +  22   \right]   }{230400}.
\end{equation}
\textcolor{black}{By expanding the above equation, we obtain}
\begin{equation}
  {\cal F}^{Shen}(k^*) =   -(k^*)^2 \left[ 1  - \frac{57 (k^*)^6}{4480} + \frac{1649 (k^*)^{8}}{737280}+O\left((k^*)^{9}\right)\right]	,
\end{equation}
which shows that Shen's scheme attains sixth order of accuracy. The operators for these schemes are plotted in Fig. \ref{fig:second_deriv}. Although Shen's scheme is sixth-order accurate, it clearly lacks damping for high-frequency modes (high wave numbers), resulting in zero damping for the highest frequency and thus the odd-even decoupling. In contrast, the sixth-order $\alpha$-damping scheme is able to resolve the higher frequency modes and provides the necessary damping. For comparison, the fourth-order $\alpha$-damping scheme of Nishikawa is also shown in Fig. \ref{fig:second_deriv}.

We emphasize that a viscous scheme has to be designed deliberately to avoid the odd-even decoupling, and meeting a desired order of accuracy is not sufficient. This is especially important for conservative schemes because high-frequency damping can be lost in the final step of the flux difference. To see this, one can express Shen's scheme in the form of a linear combination of the $m \Delta x$-spaced central differences: 
\begin{eqnarray}
 \left.
  \frac{\partial F}{\partial x} \right|_{j} 
  =  
  \sum_{m=1}^8 D_m \frac{    u_{j+m} - 2 u_j + u_{j-m}  }{  (m \Delta x)^2 }, 
\end{eqnarray}
where 
\begin{eqnarray} 
D_1 = \frac{31895}{131072},  \quad D_2 = \frac{333251}{153600},  \quad D_3 = - \frac{ 6967107}{3276800},  \quad D_4 = \frac{26711}{30720}, \\  [2ex]
D_5 = - \frac{69475}{393216},  \quad D_6 = \frac{223}{10240},  \quad D_7 =- \frac{  13769 }{ 9830400 },  \quad D_8 = \frac{3}{51200}.
\end{eqnarray}
The $m \Delta x$-spaced central differences decouple odd- and even-numbered cells when $m$ is even (as well known), thus having no high-frequency damping for the highest-frequency modes. On the other hand, they provide high-frequency damping when $m$ is odd; but we have mixed signs for $D_1$, $D_3$, $D_5$, and $D_7$, and the damping properties from these terms cancel one another, as shown in the Fourier analysis. In contrast, the new sixth-order $\alpha$-damping scheme can be written as
\begin{eqnarray}
 \left.
  \frac{\partial F}{\partial x} \right|_{j} 
  =  
  \sum_{m=1}^3 D_m \frac{    u_{j+m} - 2 u_j + u_{j-m}  }{  (m \Delta x)^2 }, 
\end{eqnarray}
where 
\begin{eqnarray} 
D_1 = \frac{3}{2},  \quad D_2 = - \frac{3}{5}, \quad D_3 =  \frac{ 1}{10}.
\end{eqnarray}
The coefficients are positive for all the odd terms: $D_1$ and $D_3$, and therefore high-frequency damping properties can never cancel each other. Indeed, this scheme provides high-frequency damping as shown in the Fourier analysis. Therefore, the explicit addition of the damping term in the numerical flux as in Equation (\ref{good_6tth-order_scheme_flux}) prevents the cancellation. Moreover, it enables the efficient construction of a high-order scheme as we have shown. On the other hand, many central-difference schemes directly approximating the diffusion term $\partial^2 u / \partial x^2 $ at a cell center often automatically have high-frequency damping. Indeed, one can write a central difference scheme as the sum of consistent and damping terms and identify the damping term for a specific scheme. See Ref.\textcolor{black}{\cite{nishikawa:AIAA2010}} for examples. 
}

\section{Results}\label{sec-4}

{\color{black} To demonstrate the importance of high-frequency damping in viscous schemes, we compare the newly derived sixth-order viscous scheme and Shen's scheme for the viscous shock-tube problem in Ref.\cite{daru2009numerical}.} In this problem, the propagation of the incident shock wave and contact discontinuity leads to the development of a thin boundary layer at the bottom wall. After its reflection at the right wall, the shock wave interacts with this boundary layer. These interactions result in a complex vortex system, separation region, and a typical lambda-shaped shock pattern. The flow contains structures of various length scales, making it an ideal test case for evaluating high-resolution schemes. The initial conditions are:
\begin{equation}\label{vst}
\begin{aligned}
\left( {\rho , u,v, p} \right) = \left\{ \begin{array}{l}
\left( {120, 0 ,0,120/\gamma } \right),  \quad 0 < x < 0.5,\\
\left( {1.2 , 0 ,0, 1.2/\gamma } \right),  \quad 0.5 \le x < 1,
\end{array} \right.
\end{aligned}
\end{equation}

with the ratio of specific heats of $\gamma = 7/5$. The domain for this test case is taken as $x \in [0,1], y \in [0,0.5]$. The flow is simulated for time $t=1$, keeping the Mach number of the shock wave at 2.37 \textcolor{black}{and constant dynamics viscosity $\mu$ = $1/500$ and $1/1000$. If the reference values are chosen as the initial speed of sound, unit density and unit length, the Reynolds number, $Re$ , will be 500 and 1000, respectively. The gas constant, $R$, is taken as 1.  The flow structures are much more complicated for these Reynolds numbers as the boundary layer separates at several points, giving rise to the development of highly complex vortex structures and interactions between vortices and shock waves. The problem is solved for these two Reynolds numbers on different grid sizes with various schemes.} The observations from the simulations are as follows:

\begin{itemize}
\item First, we show the fine grid simulations for both $Re$= 500 and 1000 on a grid size of 2000 $\times$ 1000 and 4000 $\times$ 2000, respectively, using the new sixth-order $\alpha$-damping scheme in Fig. \ref{fig:VST_fine}. Grid converged results for this test case for $Re$ = 1000 can be found in Ref. \cite{zhou2018grid} (their Fig. 6e), which is carried out on a grid resolution of 5000 $\times$ 2500. Similarly, converged results for $Re$ = 500 can be found in Ref. \cite{daru2009numerical}. The current simulations are consistent with their fine grid simulations.

\begin{figure}[H]
\centering\offinterlineskip
\subfigure[\textcolor{black}{Fine grid, $Re$ = 500.}]{\includegraphics[width=0.45\textwidth]{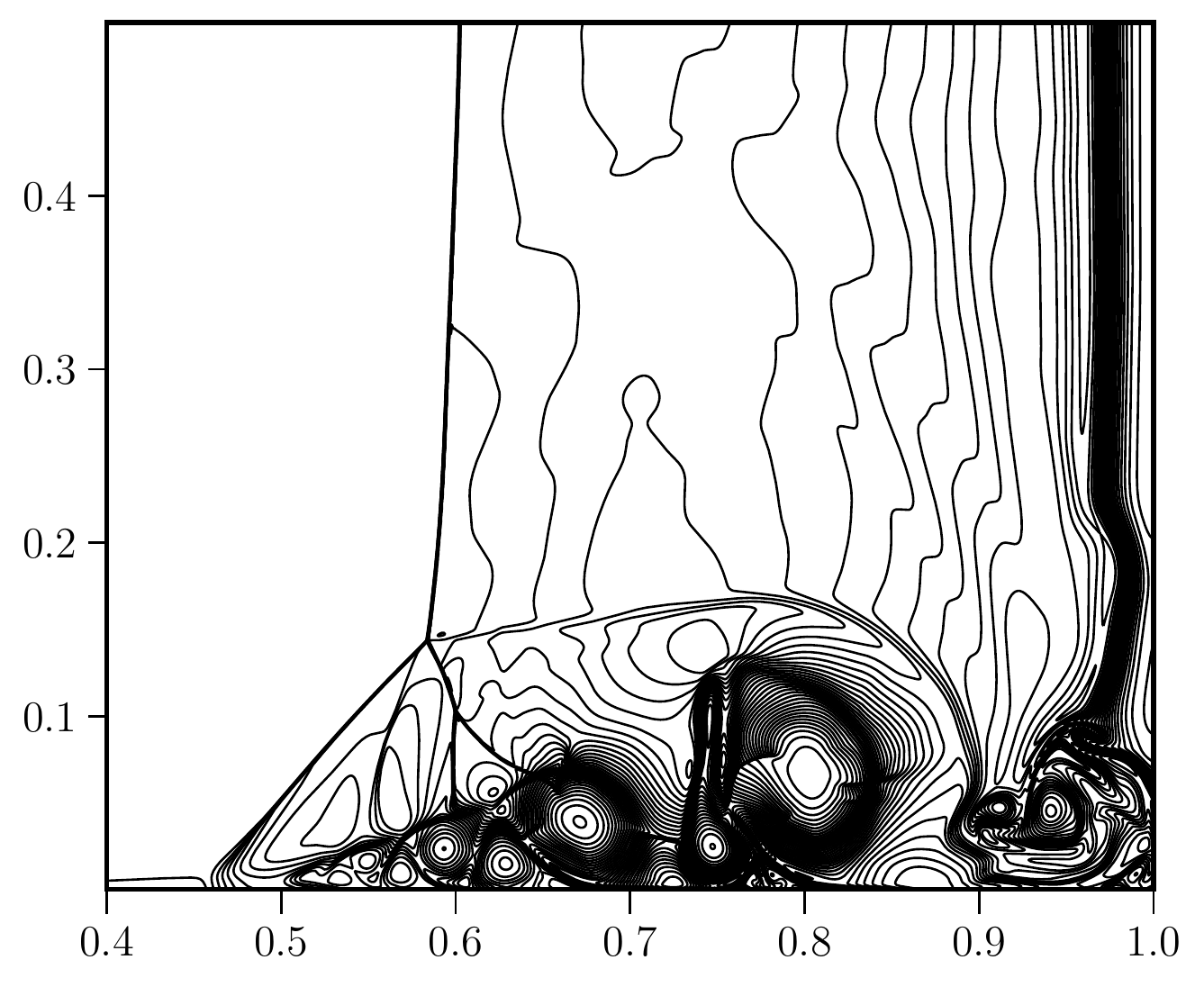}
\label{fig:500_fine}}
\subfigure[\textcolor{black}{Fine grid, $Re$ = 1000.}]{\includegraphics[width=0.45\textwidth]{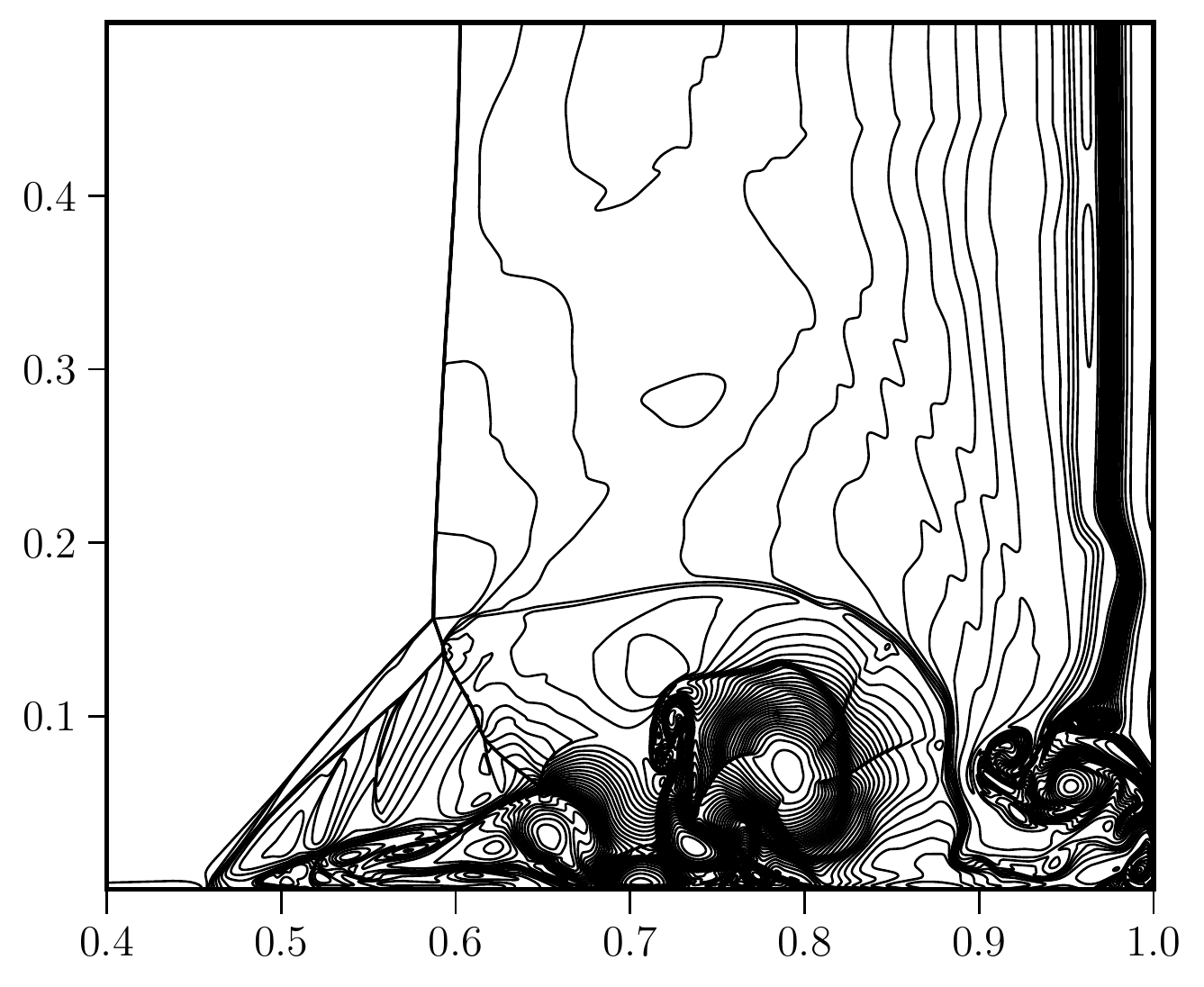}
\label{fig:1000_fine}}
\caption{\textcolor{black}{Density contours of the fine grid simulations for the viscous shock tube test case using the proposed viscous scheme. Fig. (a) $Re$ = 500 on a grid size of 2000 $\times$ 1000 and Fig. (b) $Re$ = 1000 on a grid size of 4000 $\times$ 2000.}}
\label{fig:VST_fine}
\end{figure}
\item Next, we show the effect of \textcolor{black}{viscous flux computations using Shen's scheme}, Equation (\ref{bad_6th-order_scheme_flux}), and the $\alpha$-damping approach, Equation (\ref{good_6tth-order_scheme_flux}), on the solution quality for $Re$ = 500. Fig. \ref{fig_500_t1} shows all the simulations computed on a grid size of 640 $\times$ 320.  From Fig. \ref{fig:500_shen} we can observe that the solution obtained by \textcolor{black}{Shen's scheme} shows a significant difference from the converged results in Fig. \ref{fig:500_fine}, with the primary vortex being completely distorted. On the other hand, the solution obtained by the $\alpha$-damping approach is similar to the converged result despite using a coarser grid. This clearly shows the benefits of the $\alpha$-damping approach.

\item \textcolor{black}{Figs. \ref{fig_500_t1} are at time $t$ = 1 and indicate no odd-even decoupling or high-frequency errors but if we look at the intermediate times $t$ = 0.45, 0.50, 0.54 and 0.65 shown in Fig.\ref{fig:500_shen_nodamp_45}-\ref{fig:500_shen_nodamp_65}, one can see clear saw-tooth like oscillations due to the lack of high-frequency damping. Such oscillations are not visible in Fig.\ref{fig:500_4e6_nodamp_45}-\ref{fig:500_4e6_nodamp_65} which includes the high-frequency damping \textcolor{black}{through the $\alpha$-damping approach}.}

\item \textcolor{black}{ An important observation here is that even though both the schemes are sixth-order accurate, the $\alpha$-damping approach gave better results due to the difference in spectral properties.} These results indicate that the spectral properties of the \textcolor{black}{viscous fluxes} do affect the solution quality of viscous flows, and the viscous flux discretisation is equally as important as inviscid fluxes and should be carefully evaluated.

\begin{figure}[H]
\centering\offinterlineskip
\subfigure[\textcolor{black}{$Re$ = 500, Viscous scheme of Shen}]{\includegraphics[width=0.45\textwidth]{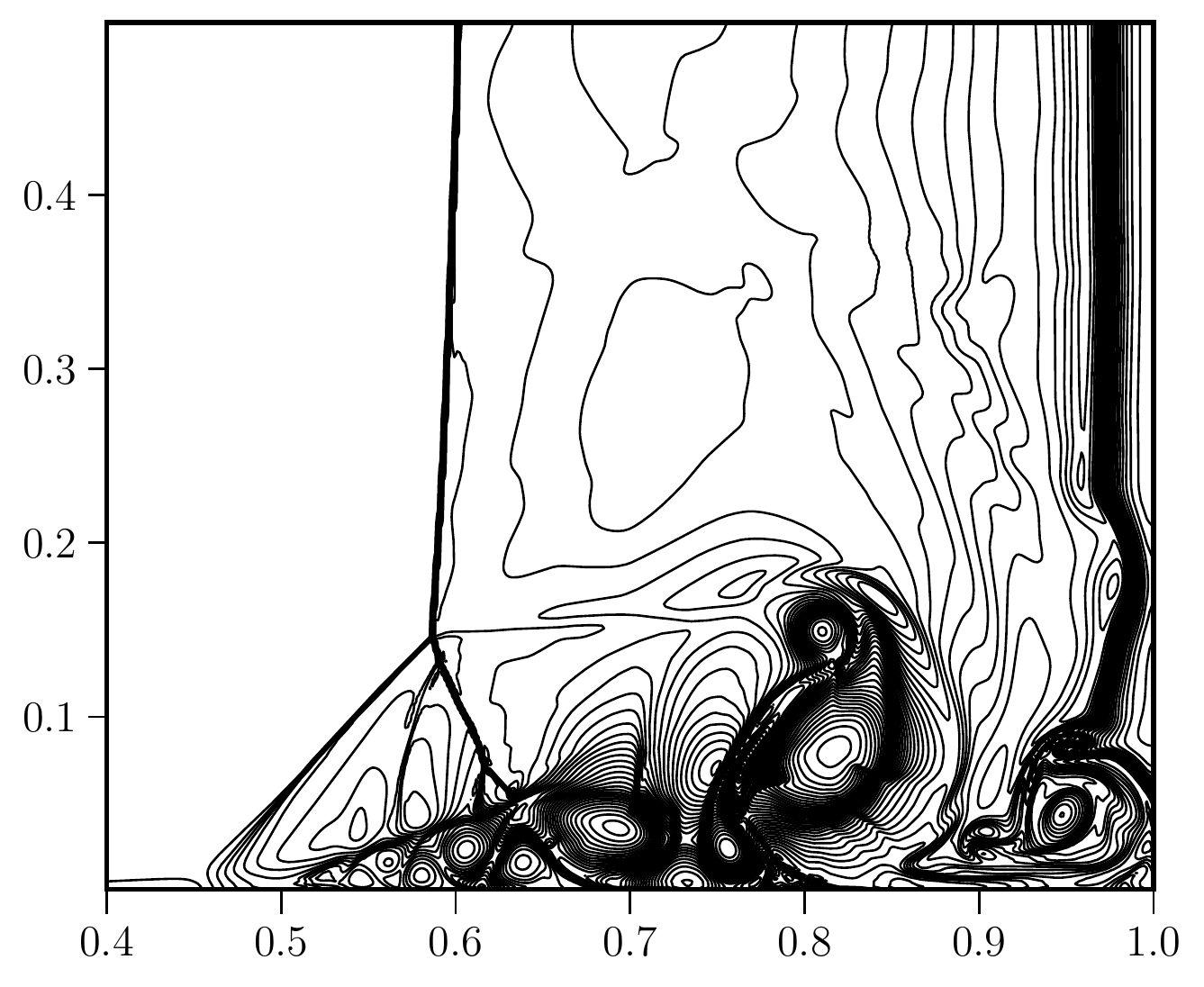}
\label{fig:500_shen}}
\subfigure[\textcolor{black}{$Re$ = 500, $\alpha$-damping}]{\includegraphics[width=0.45\textwidth]{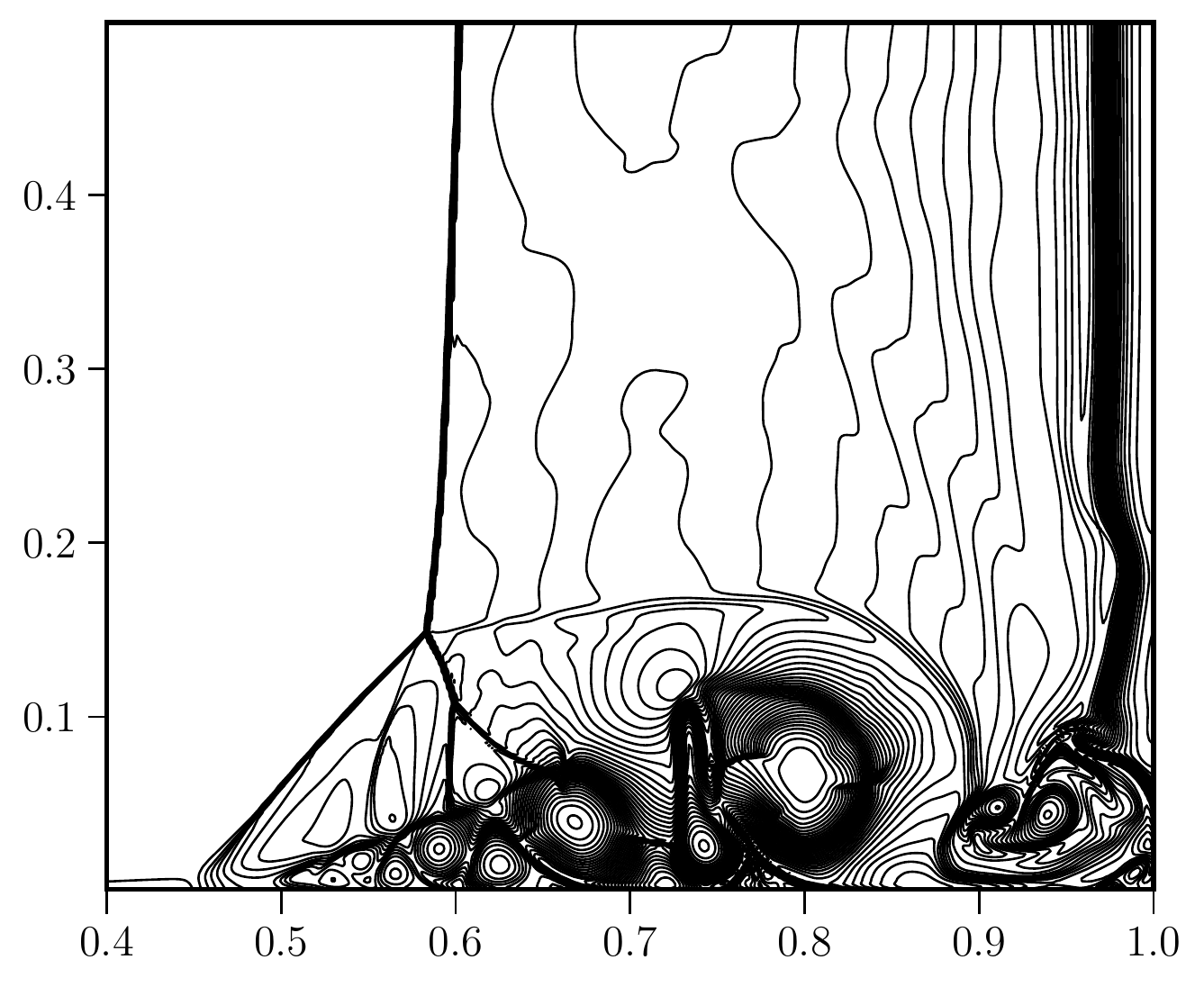}
\label{fig:500_4e6}}
\caption{\textcolor{black}{Density contours for $Re=500$ on a grid size of 640 $\times$ 320 using the viscous scheme of Shen and the $\alpha$-damping approach.}}
\label{fig_500_t1}
\end{figure}

\begin{figure}[H]
\begin{onehalfspacing}
\centering\offinterlineskip
\subfigure[\textcolor{black}{$t$ = 0.45, Shen's scheme}]{\includegraphics[width=0.15\textheight]{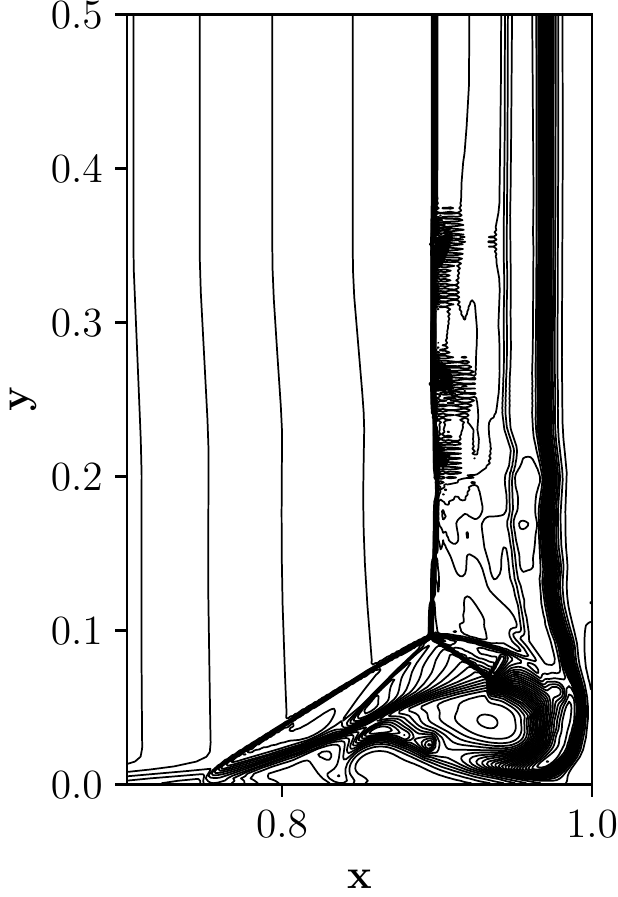}
\label{fig:500_shen_nodamp_45}}
\subfigure[\textcolor{black}{$t$ = 0.50, Shen's scheme}]{\includegraphics[width=0.15\textheight]{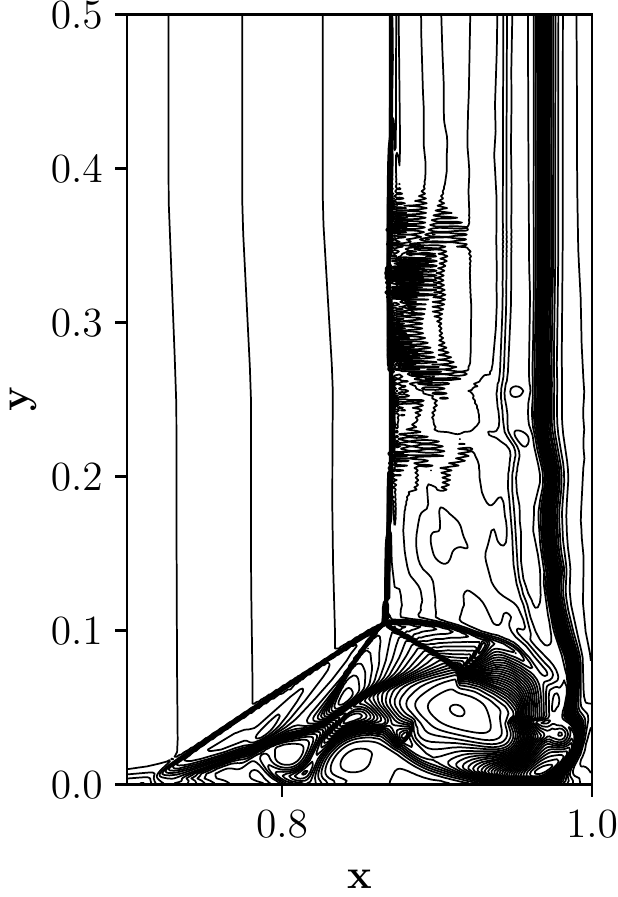}
\label{fig:500_shen_nodamp_50}}
\subfigure[\textcolor{black}{$t$ = 0.54, Shen's scheme}]{\includegraphics[width=0.15\textheight]{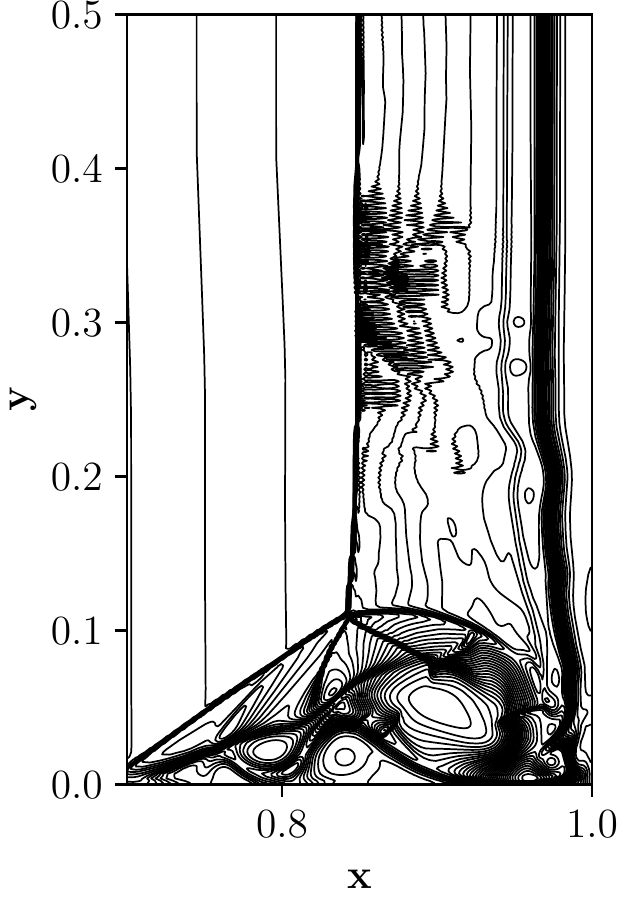}
\label{fig:500_shen_nodamp_54}}
\subfigure[\textcolor{black}{$t$ = 0.65, Shen's scheme}]{\includegraphics[width=0.15\textheight]{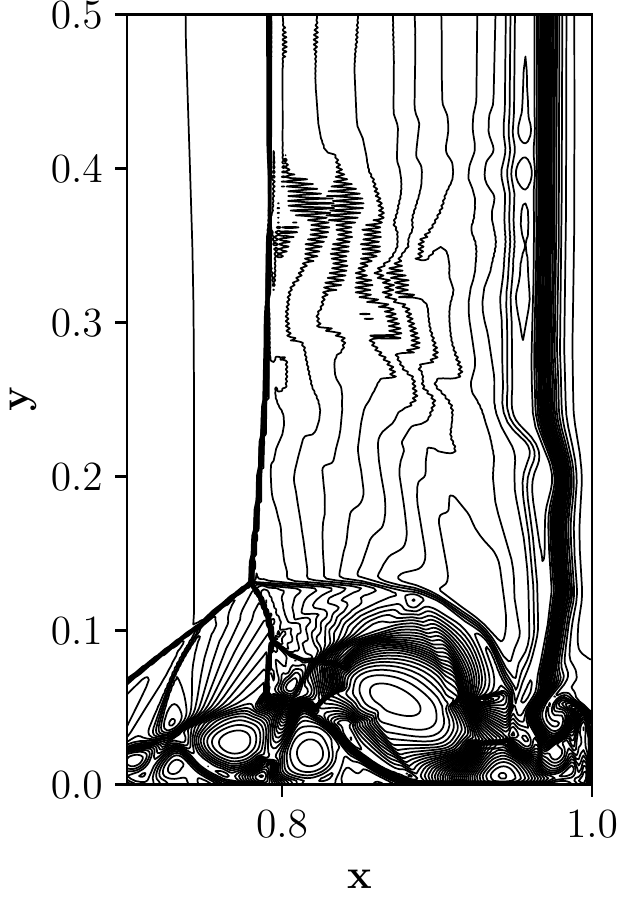}
\label{fig:500_shen_nodamp_65}}

\subfigure[\textcolor{black}{$t$ = 0.45, $\alpha$-damping}]{\includegraphics[width=0.15\textheight]{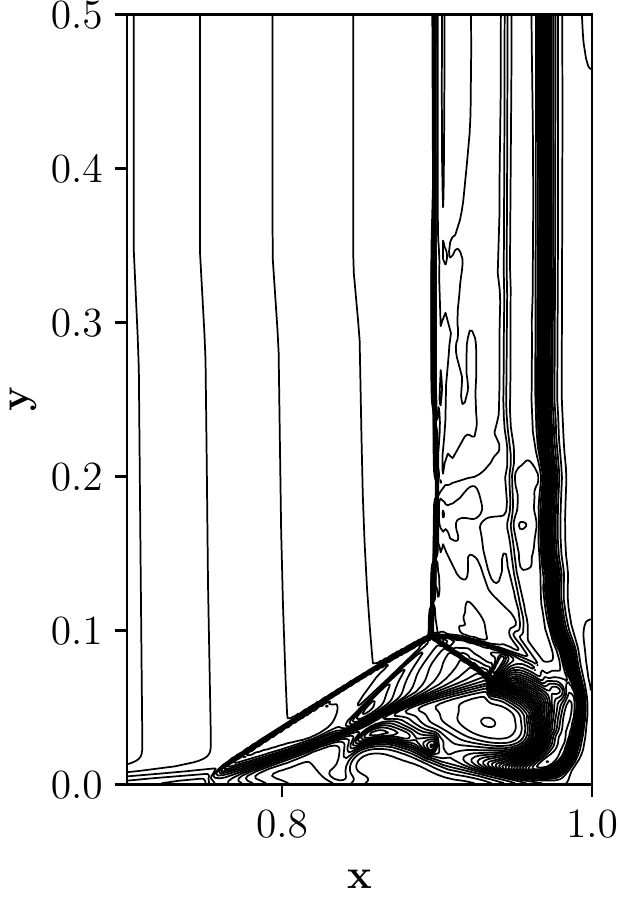}
\label{fig:500_4e6_nodamp_45}}
\subfigure[\textcolor{black}{$t$ = 0.50, $\alpha$-damping}]{\includegraphics[width=0.15\textheight]{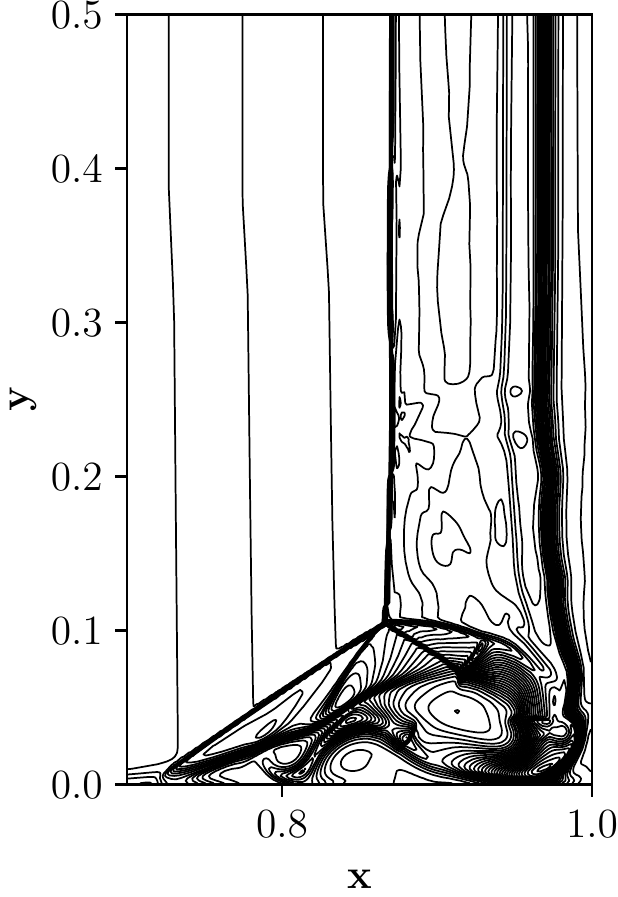}
\label{fig:500_4e6_nodamp_50}}
\subfigure[\textcolor{black}{$t$ = 0.54, $\alpha$-damping}]{\includegraphics[width=0.15\textheight]{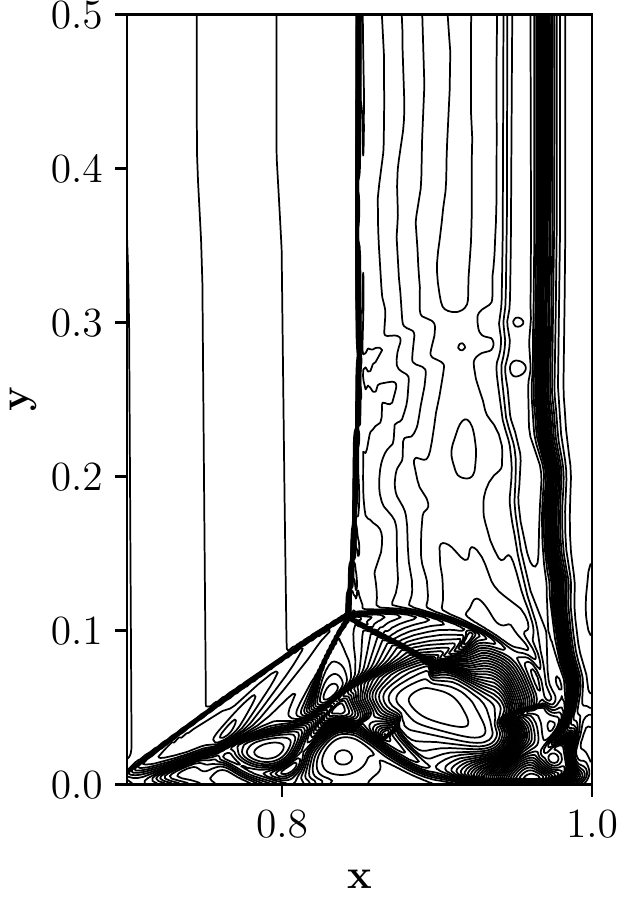}
\label{fig:500_4e6_nodamp_54}}
\subfigure[\textcolor{black}{$t$ = 0.65, $\alpha$-damping}]{\includegraphics[width=0.15\textheight]{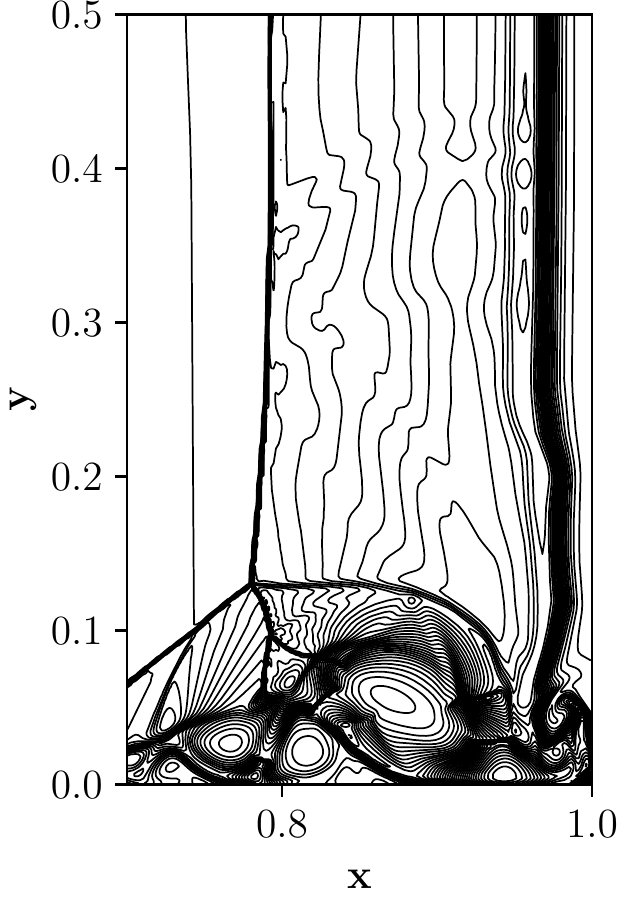}
\label{fig:500_4e6_nodamp_65}}
\caption{\textcolor{black}{Density contours for intermediate times for $Re=500$ on a grid size of 640 $\times$ 320 using the viscous scheme of Shen, top row, and the $\alpha$-damping approach, bottom row. These figure are drawn with 38 density contours.}}
\label{fig_damp_500}
\end{onehalfspacing}
\end{figure}

\item Similar observations can be made for $Re$ = 1000, where simulations are carried out on a grid size of 1280 $\times$ 640. For time $t$ =1 the solution obtained by the $\alpha$-damping approach, Fig. \ref{fig:1000_4e6}, is closer to the converged results, shown in Fig. \ref{fig:1000_fine}, than Shen's scheme Fig. \ref{fig:1000_shen}. The lambda-shaped shocks are also captured accurately along with the small scale features. The rotating vortices on the lower right corner also fit very well with the converged results. Similar to the observations made earlier for $Re$ = 500 using Shen's scheme, the intermediate solutions for the $Re$ =1000 also show oscillations as shown in Figs. \ref{fig:1000_shen_nodamp_45}-\ref{fig:1000_shen_nodamp_65}. On the other hand, the $\alpha$-damping approach shows no oscillations at all.

\begin{figure}[H]
\centering\offinterlineskip
\subfigure[\textcolor{black}{$Re$ = 1000, Shen's scheme}]{\includegraphics[width=0.45\textwidth]{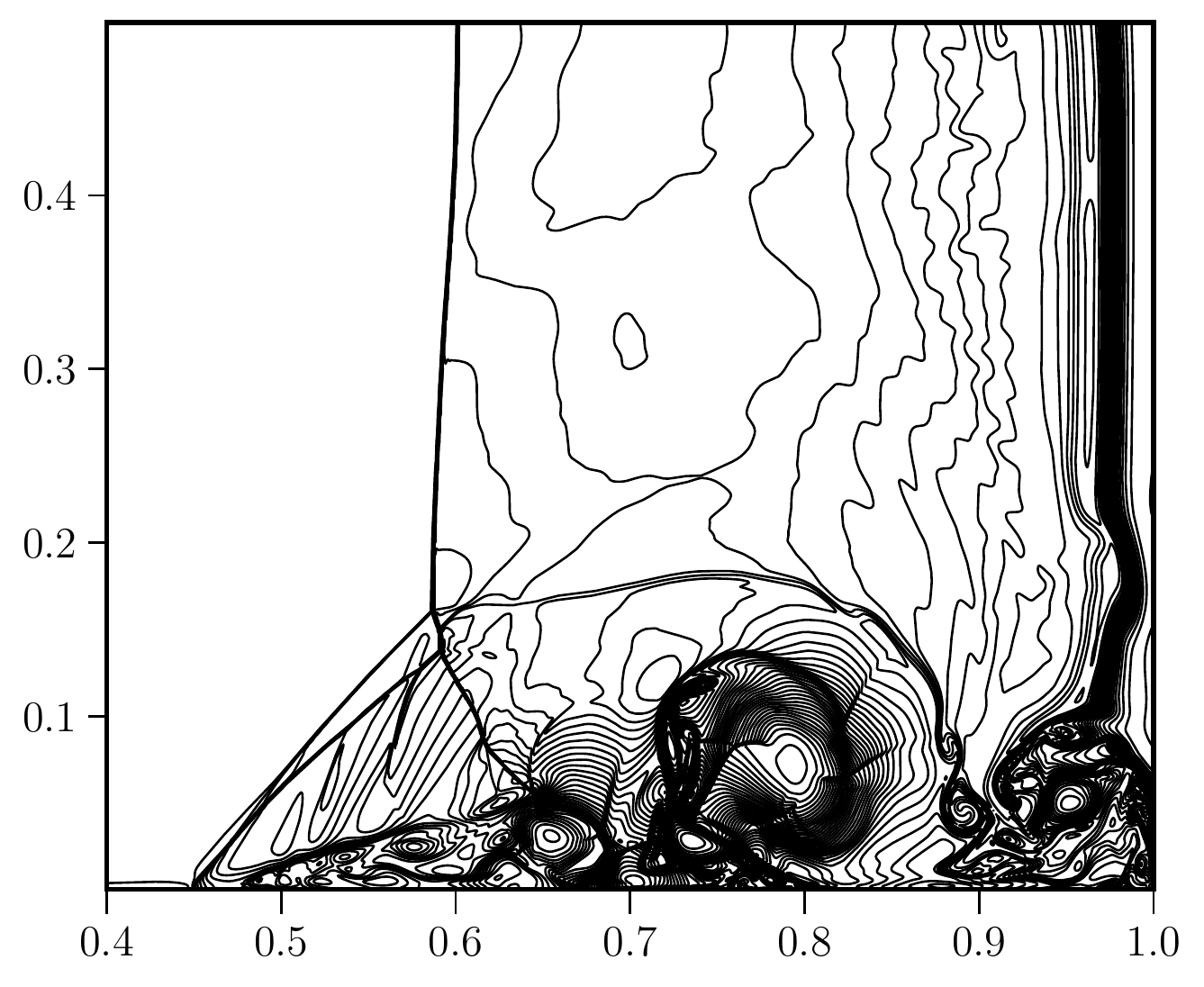}
\label{fig:1000_shen}}
\subfigure[\textcolor{black}{$Re$ = 1000, $\alpha$-damping}]{\includegraphics[width=0.45\textwidth]{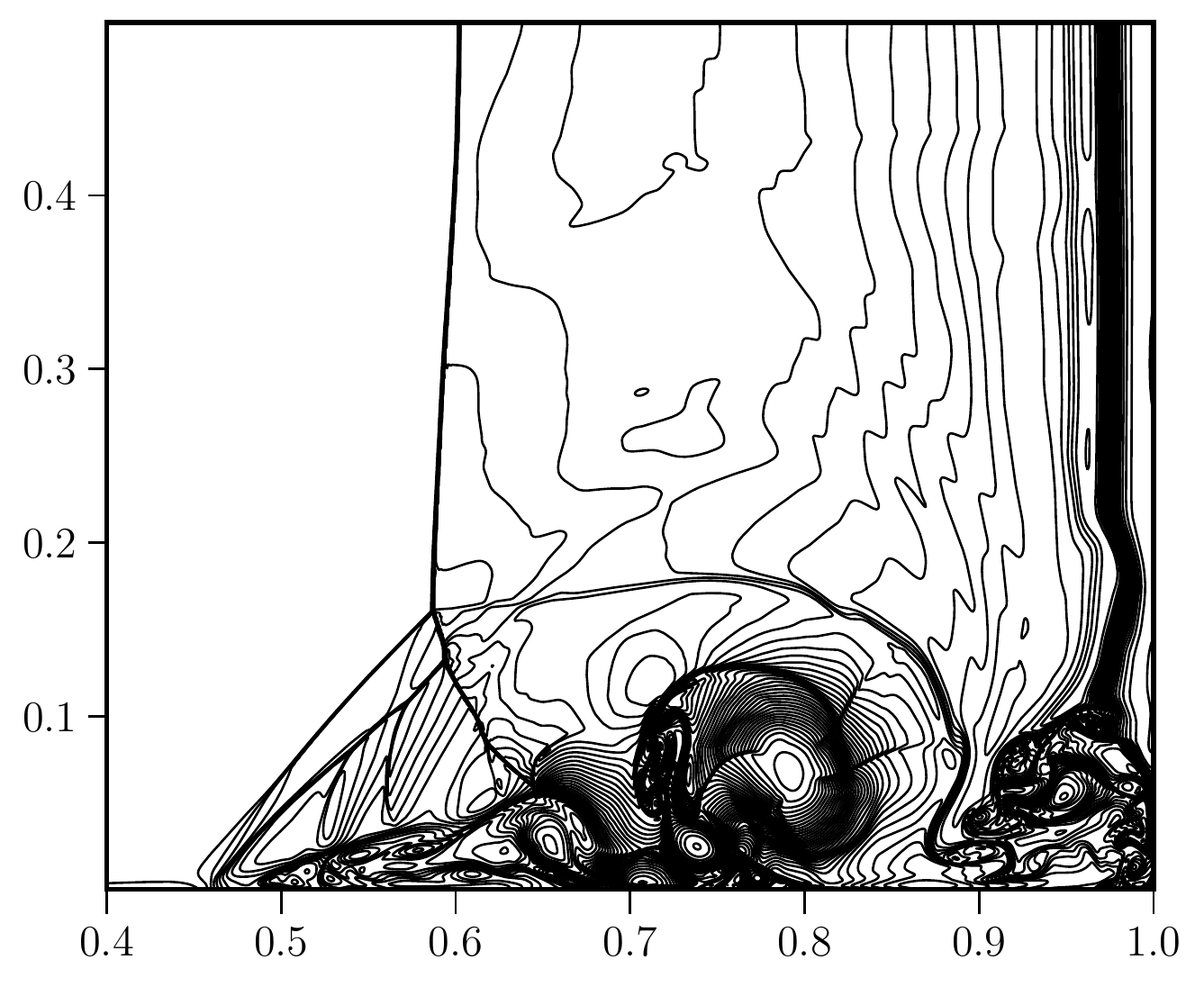}
\label{fig:1000_4e6}}
\caption{\textcolor{black}{Density contours for $Re=1000$ on a grid size of 1280 $\times$ 640 using Shen's scheme and the $\alpha$-damping approach.}}
\label{fig_1000_t1}
\end{figure}

\begin{figure}[H]
\begin{onehalfspacing}
\centering\offinterlineskip
\subfigure[\textcolor{black}{$t$ = 0.45, Shen's scheme}]{\includegraphics[width=0.15\textheight]{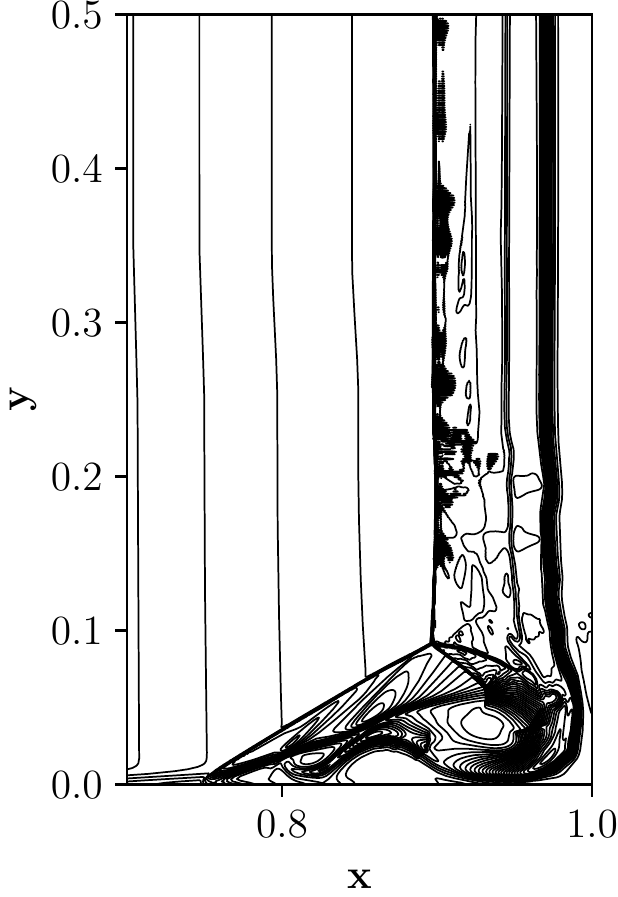}
\label{fig:1000_shen_nodamp_45}}
\subfigure[\textcolor{black}{$t$ = 0.50, Shen's scheme}]{\includegraphics[width=0.15\textheight]{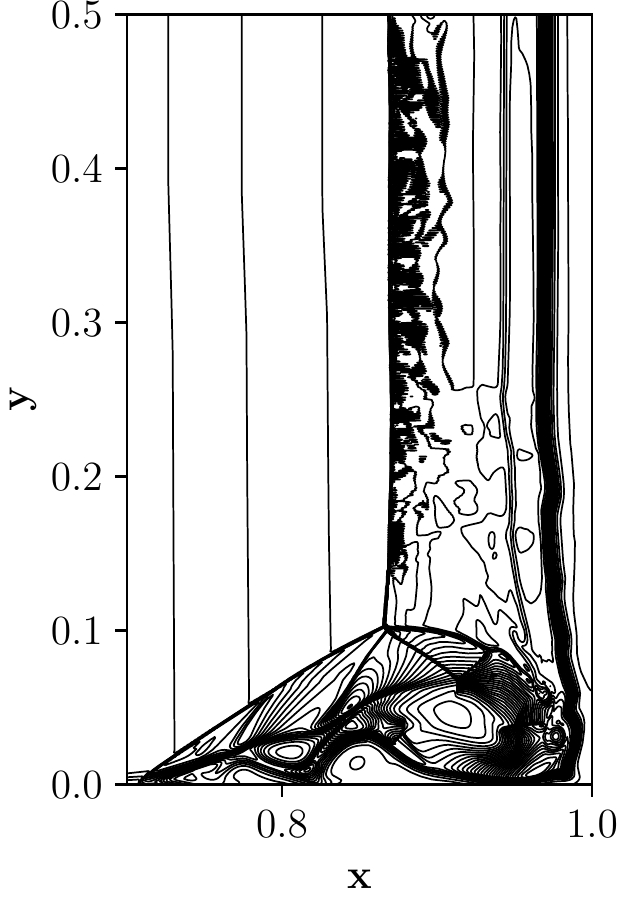}
\label{fig:1000_shen_nodamp_50}}
\subfigure[\textcolor{black}{$t$ = 0.54, Shen's scheme}]{\includegraphics[width=0.15\textheight]{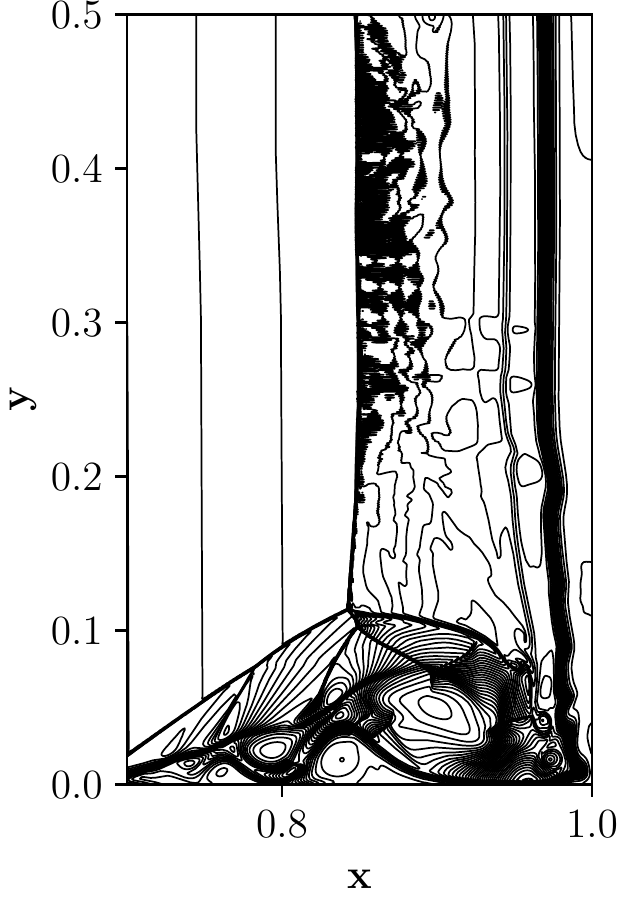}
\label{fig:1000_shen_nodamp_54}}
\subfigure[\textcolor{black}{$t$ = 0.65, Shen's scheme}]{\includegraphics[width=0.15\textheight]{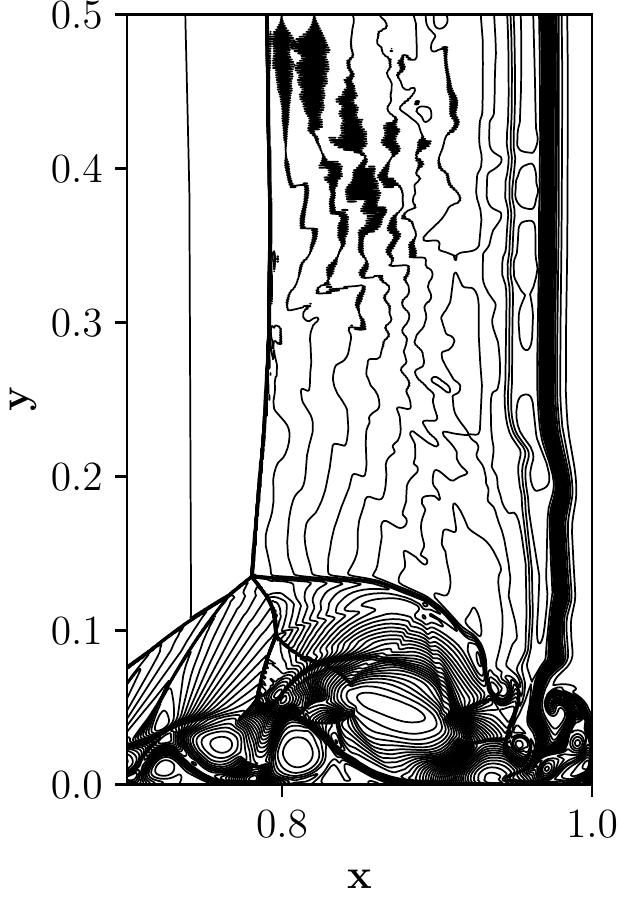}
\label{fig:1000_shen_nodamp_65}}

\subfigure[\textcolor{black}{$t$ = 0.45, $\alpha$-damping}]{\includegraphics[width=0.15\textheight]{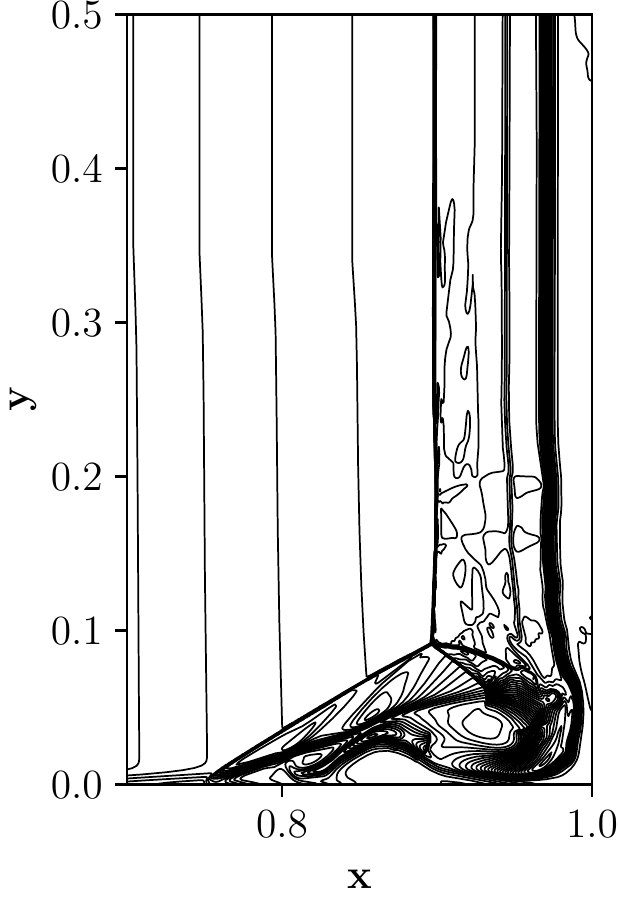}
\label{fig:1000_4e6_nodamp_45}}
\subfigure[\textcolor{black}{$t$ = 0.50, $\alpha$-damping}]{\includegraphics[width=0.15\textheight]{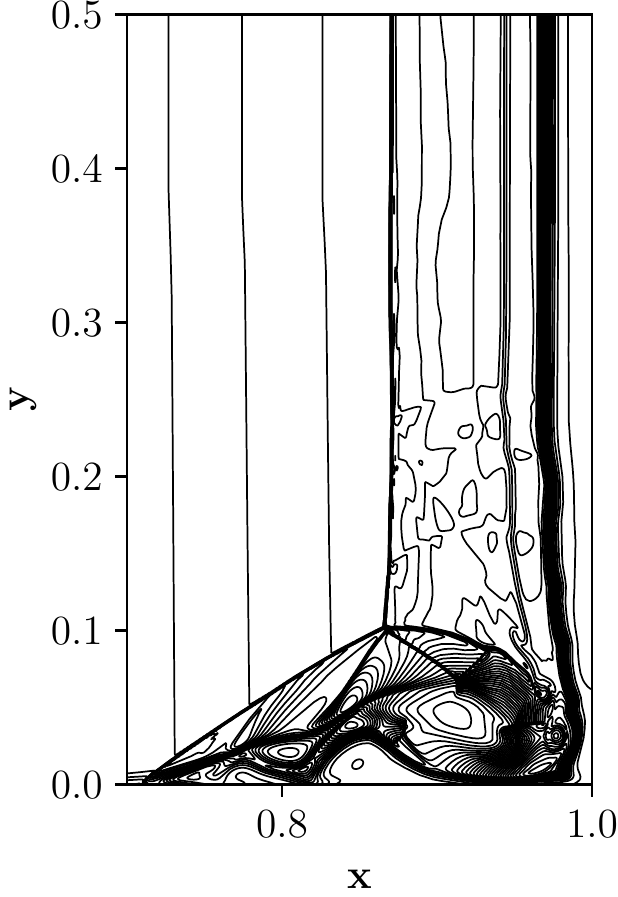}
\label{fig:1000_4e6_nodamp_50}}
\subfigure[\textcolor{black}{$t$ = 0.54, $\alpha$-damping}]{\includegraphics[width=0.15\textheight]{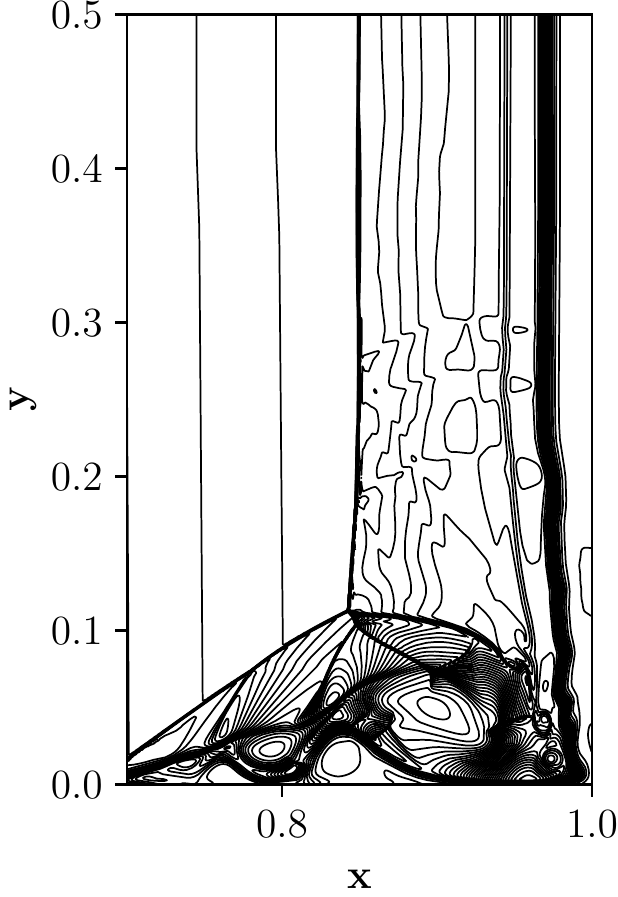}
\label{fig:1000_4e6_nodamp_54}}
\subfigure[\textcolor{black}{$t$ = 0.65, $\alpha$-damping}]{\includegraphics[width=0.15\textheight]{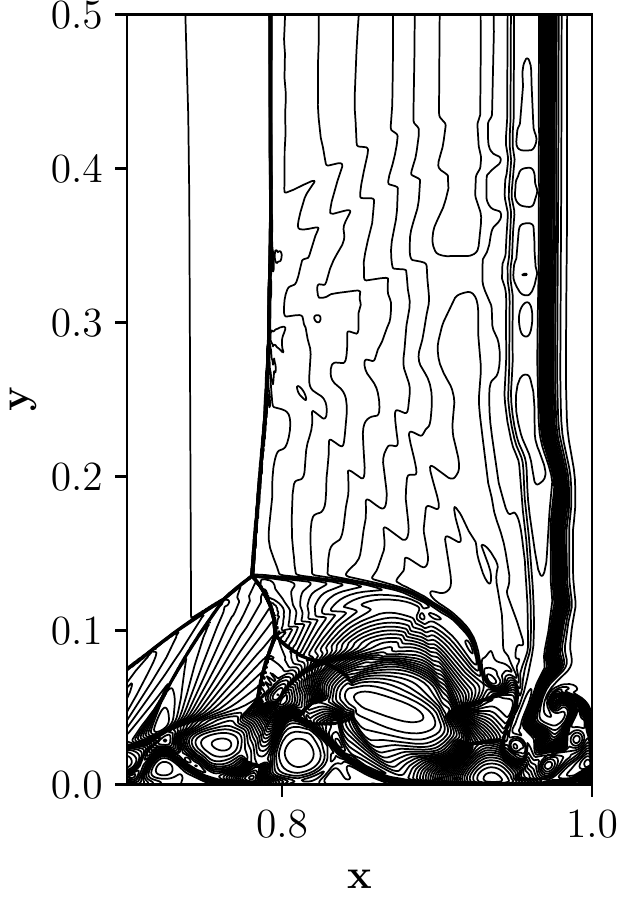}
\label{fig:1000_4e6_nodamp_65}}
\caption{\textcolor{black}{Density contours for intermediate times for $Re=1000$ on a grid size of 1280 $\times$ 640 using Shen's scheme, top row, and the $\alpha$-damping approach, bottom row. These figure are drawn with 38 density contours.}}
\label{fig_damp_1000}
\end{onehalfspacing}
\end{figure}
\end{itemize}
  
\section{Conclusions}\label{sec-5} 

\textcolor{black}{The results presented in this note demonstrate the importance of high-frequency damping in the discretization of the viscous fluxes. As we have shown, taking Shen's scheme as an example, conservative discretizations can easily lack the high-frequency damping property (thus leading to odd-even decoupling) if the diffusive flux at a face is designed only to achieve a desired order of accuracy and not with a damping term explicitly incorporated. Although a damping term can be introduced to such schemes, we have shown that there is a simpler approach, whereby extending the fourth-order $\alpha$-damping scheme in Ref.\textcolor{black}{\cite{nishikawa:AIAA2010}}, we can construct a sixth-order scheme by carefully combining lower-order gradients and a damping term. While formally sixth-order accurate, a scheme lacking a high-frequency damping property has been shown to generate high-frequency oscillations, implying the potential to provide errors in viscous simulations, for example the oscillations can potentially induce false transition to turbulence. On the other hand, the efficient sixth-order $\alpha$-damping scheme derived in this note has been shown to suppress such oscillations and generate results in very good agreement with reference solutions. These results indicate that a properly-designed viscous discretization is essential to accurate simulations and can be just as impactful on the final solution as the inviscid flux discretization.}

\section*{Acknowledgements}
\textcolor{black}{This work has received funding from the European Union’s Horizon 2020 Research And Innovation Programme under grant agreement No. 815278. The authors gratefully acknowledge the financial support provided by the European Union’s Horizon 2020 Research And Innovation Programme under grant agreement No. 815278. They also thank the support given by the Technion-Israel Institute of Technology-Haifa, Israel.} 


\bibliographystyle{elsarticle-num}
\bibliography{mig_ref}

\end{document}